\documentclass[smallextended,referee,envcountsect]{svjour3}
\UseRawInputEncoding
\smartqed
\usepackage{graphicx}
\journalname{JOTA}

\def\ra{\rightarrow}

\def\a{\alpha}

\def\s{\sigma}

\def\n{\nabla}

\def\r1{\mathbb{R}}
\def\mb{\mathbb}

\def\mc{\mathcal}
\def\mb{\mathbb}

\def\be{\begin{array}}
\def\en{\end{array}}
\def\befc{ \begin{equation}\left\{\begin{array}{llll}}
\def\enfc{ \end{array} \right. \end{equation}}
\def\beeq{\begin{equation}}
\def\eneq{\end{equation}}
\def\beit{\begin{itemize}}
\def\enit{\end{itemize}}
\def\bepm{\begin{pmatrix}}
\def\enpm{\end{pmatrix}}

\def\htheta{\mbox{\boldmath$\theta$}}   
\def\hx{\mbox{\boldmath$x$}}   
\def\hlambda{\mbox{\boldmath$\lambda$}}   
\def\hnu{\mbox{\boldmath$\nu$}}   
\def\hL{\mbox{\boldmath$\mathcal{L}$}}   
\def\hw{\mbox{\boldmath$w$}}   
\def\hh{\mbox{\boldmath$h$}}   

\usepackage{mathrsfs}
\usepackage[compress]{cite}
\usepackage{amsmath,amssymb,amsfonts}
\usepackage{algorithmic}
\usepackage{graphicx}
\usepackage{subfigure}
\usepackage{textcomp}
\usepackage{empheq}
\usepackage{mathrsfs}
\usepackage{color}

\allowdisplaybreaks[4]   

\newtheorem{assumption}{Assumption}{}

\usepackage[dvips,colorlinks,dvipdfm]{hyperref}
\hypersetup{colorlinks,
linkcolor=red,
citecolor=red,
urlcolor=blue}

\makeatletter

\newcommand{\Rmnum}[1]{\expandafter\@slowromancap\romannumeral #1@}
\makeatother

\newenvironment{assbis}
  {\addtocounter{assumption}{-2}%
   \begin{assumption}}
  {\end{assumption}}

\makeatletter
\def\rightharpoonfill@{\arrowfill@\relbar\relbar\rightharpoonup}
\DeclareRobustCommand{\overrightharpoon}{\mathpalette{\overarrow@\rightharpoonfill@}}
\makeatother

\begin{document}

\title{A Multi-Scale Method for Distributed Convex Optimization with Constraints}

\subtitle{}

\author{Wei Ni \and Xiaoli Wang }

\institute{Wei Ni, Corresponding author \at
             School  of Sciences, Nanchang University \\
              Nanchang 330031, China \\
              niw@amss.ac.cn
           \and
           Xiaoli Wang\at
              School of Information Science and Engineering\\
              Harbin Institute of Technology at Weihai \\
              Weihai 264209, China \\
              xiaoliwang@amss.ac.cn
}

\date{Received: date / Accepted: date}

\maketitle

\begin{abstract}
This paper proposes a multi-scale method to design a continuous-time
distributed algorithm for constrained convex optimization problems by using multi-agents with Markov switched network dynamics and noisy inter-agent communications. Unlike most previous work which mainly puts emphasis on dealing with fixed network topology, this paper tackles the challenging problem of investigating the joint effects of stochastic networks and the inter-agent communication noises on the distributed optimization dynamics, which has not been systemically studied in the past literature. Also, in sharp contrast to previous work in constrained optimization, we depart from the use of projected gradient flow which is non-smooth and hard to analyze; instead, we  design a smooth optimization dynamics which leads to easier convergence analysis and more efficient numerical simulations. Moreover, the multi-scale method presented in this paper generalizes previously known distributed convex optimization algorithms from the fixed network topology to the switching case and the stochastic averaging obtained in this paper is a generalization of the existing deterministic averaging.
\end{abstract}


\section{Introduction}

The research of convex optimization by using multi-agent systems is a hot topic in recent decades. This problem typically  takes the form of  minimizing the sum of functions and  is usually divided into subtasks of local optimizations, where each local optimization subtask is  executed by one agent and  the cooperation among these agents makes these local algorithms compute the optimal solution in a consensus way.
Informally, convex optimization algorithms constructed  in this way are usually termed as distributed convex optimization (DCO), which models a broad array of engineering and economic scenarios and finds numerous applications in diverse areas such as operations research, network flow optimization, control systems and signal processing; see \cite{nedic2009,nedic2010,nedic2015,boyd2011,feijer2010, duchi2012, ram2010}  and references therein for more details.

Usually, DCO takes advantages of the consensus algorithms  in multi-agent systems and gradient algorithms in convex optimization.
The idea of combining them was proposed early in 1980s by Tsitsiklis et al. in \cite{tsitsiklis1986} and re-examined recently in the context of DCO  in \cite{nedic2009,nedic2010,boyd2011,duchi2012}.
Most DCO algorithms in the earlier development were discrete-time, with the distributed gradient descent strategy \cite{nedic2009, nedic2010, ram2010} being the most popular.
Various extensions of the  distributed gradient descent  were then proposed, such as push-sum based approach \cite{nedic2015}, incremental gradient procedure \cite{ram2009}, proximal method \cite{parikh2013}, fast distributed gradient strategy \cite{jakovetic2014,qu2017},
and non-smooth analysis based technique \cite{zeng2017}, just to name a few.  Additionally, using local gradients is rather slow and the community has moved towards using gradient estimation \cite{xin2018,qu2017} and stochastic gradient \cite{ram2010, yuan2016}.    Further research directions were then followed by taking  optimization constraints into consideration. Generally, DCO with constraints has proceeded along two research lines, namely projected gradient strategy  and primal-dual scheme. The projected gradient strategy designed the optimization algorithms by projecting the gradient into the constraint set and extended it by including a consensus term \cite{ram2010,duchi2012,lou2016}.  The primal-dual scheme introduced equality and inequality multipliers and designed for them extra dual  dynamics  \cite{charalamous2014,yi2015,zhu2011} in which a projection onto positive quadrant is  usually included \cite{feijer2010,yi2015,charalamous2014,yamashita2020} so that the inequality-multiplier stays positive. Therefore, both research lines  dealing with optimization constraints above  are projection-dependent.

While projection-dependent algorithms were widely used in DCO, they  require the optimization constraints to have a relatively simple form so that  the projections can be computed analytically.  To overcome this difficulty,   this paper pursues a new method to design a novel DCO algorithm  by avoiding projection. Our method is built on the primal-dual setup by introducing Lagrange multipliers. We modify the classical projection-based dynamics for the inequality-multiplier by utilizing the technique of mirror descent \cite{nemirovsky1983, raginsky2012} and design a projection-free multiplier-dynamics which is smooth.   In conclusion, compared with most existing constrained DCO algorithms, our method avoids projection and thus reduces the difficulties of convergence analysis and iterative computation.

Aside from the difficulty  associated with optimization constraints, the second challenge in DCO problem  ties with stochastic  networks and inter-communication noises. This challenge, together with optimization constraints,  jointly make the DCO problem difficult to analyze, and therefore relatively few results  were reported.  Nedic \cite{nedic2015} considered DCO over deterministic and  uniformly strongly connected time-varying  networks without considering communication noises and optimization constraints; furthermore their algorithm needed a strong  requirement that  each node knows its out-degree at all times.  The DCO problem with optimization constraints over time-varying graphs was investigated in \cite{xie2018} by using the epigraph form, but communication noises were not considered there.
The work in  \cite{lobel2011} also investigated the consensus-based DCO algorithm by using a random graph model where the communication link availability is described by a stochastic process, but leaving challenge issues of optimization constraints and communication noise untouched.

To tackle the above-mentioned difficulties and to contribute to the existing literature, this paper proposes a multi-scale method for constrained DCO problem over Markov switching networks under noisy communications. Unlike most existing DCO algorithms which are discrete-time, we study this problem in a continuous-time  framework  because the classical tools of Ito formula, backward Kolmogorov equation and ergodic theory in stochastic analysis can be used and the elegant Lyapunov argument in optimization theory \cite{feijer2010} can be invoked.
Recently, we established in  \cite{ni2016b} a new technique of stochastic averaging (SA) for unconstraint DCO, where the idea of averaging was perviously explored by us to handle the switching networks of the multi-agent systems, with the deterministic  version being presented in \cite{ni2013,ni2012} and the stochastic version in  \cite{ni2016a}.
Compared  with \cite{ni2016b}, the present work considers a more general case of constrained optimization which is more challenge.

Although the   SA viewpoint  for multi-agent systems has been indicated in \cite{ni2016a},  its theoretical clarification and design details were not provided there. In this paper, we generalize the SA principle in \cite{ni2016a}  to the DCO problem and  propose a multi-scale based design procedure for the SA.
We begin with the intuition behind our approach.
Our multi-scale method borrows the idea of the slave principle in Synergetics \cite{haken1982} which was initially proposed by German physician Haken in 1970s. According to this principle, the system variables are classified into fast and slow ones, where the slow variables dominate the system evolution and characterize the  ordering degree of the system. Therefore,  it is necessary to eliminate the fast variables and obtain an equation  for the slow variables only. This equation is called the principle equation and it can be viewed as an approximate description of the system.  The method of eliminating the fast variables in physics is termed as Born-Oppenheimer approximation. In this paper, the continuous-time Markov chain  characterizing the time-varying networks is regarded as the fast variable and, in contrast, the states of the optimization multi-agent system are considered as slow variables.
To distinguish the fast and slow variables,  two time-scales for them are introduced.
In this paper, we propose a concrete scheme to eliminate the fast variable by resorting to the tool of multiscale analysis introduced by Pavliotis \cite{pavliotis} and obtain an averaged SDE which acts as an approximation to the original  switching stochastic differential equation (SDE). In this sense, the effect of  network switching on optimization dynamics is eliminated, and thus the DCO under switching networks is, in fact,  reduced to that under fixed case.

We mention two benefits of our multi-scale method. The first benefit    in comparison with those dealing with random networks  (see e. g. \cite{nedic2009},\cite{lobel2011}) lies in its generalizability.
The widely used DCO algorithms over random networks in \cite{nedic2009},\cite{lobel2011} built their analysis on the product theory of stochastic matrices (c.f. \cite{wolfowitz1963}): the products of stochastic matrices converges  to a rank one matrix (see Lemmas 4-7 in \cite{lobel2011}). This theory was used to analyze  the convergence of their optimization algorithm which is  driven by a chain of stochastic matrices and an inhomogeneous gradient term (see Section V in \cite{lobel2011}).
However, due to technicalities involved, this method is hard to generalize to include optimization constraints or communication noises since otherwise the resulting matrices are not stochastic matrices. Our method does not have this limitation, instead it can treat the optimization constraints, communication noises and stochastic networks into a unified framework. As the second benefit, our method reduces the optimization algorithm in stochastic networks  to that in fixed network (see Theorem \ref{theorem4}) and establishes an approximation relationship between the two algorithms (see Theorem \ref{theorem5}). Therefore, the SA method in this paper can help generalize existing DCO algorithms from fixed network to stochastic networks.

To sum up, the contributions of our paper are as follows. Firstly,  this paper proposes a novel method of SA for the design and analysis of DCO problem, and addresses in a unified framework the challenge issues  of the optimization constraints, communication noises and stochastic  networks. Secondly, the SA method can help generalize some DCO algorithms from fixed network to switching case since it converts the latter to the former and establishes an approximation relationship between them.
Thirdly, the DCO algorithm in this paper is projection-free and thus has the advantages of removing the difficulty of computing projection  and  rending  the resulting optimization dynamics to be  smooth so that algorithm analysis and simulation become relatively easy.
Lastly,  the multi-scale method used in this paper has the ability to generalize the averaging method from the deterministic case \cite{ni2013,ni2012} to the stochastic case in the present form, and it can also provide a theoretical justification for our original vision of SA for multi-agent systems \cite{ni2016a}.

\section{Preliminaries}

{\bf A. Notations.}
For a vector $a$, its $i$-th component is denoted by $[a]_i$ or $a_i$. By $a \prec 0$ ($a\preccurlyeq 0$) we mean that each entry of $a$ is less than (less than or equal to) zero.
Letting \scalebox{0.9}[1]{$a=(a_1, \cdots, a_n)^T\in \mathbb{R}^n$} and \scalebox{0.9}[1]{$b=(b_1, \cdots, b_n)^T\in \mathbb{R}^n$}, we define $a\odot\ b=(a_1b_1, \cdots, a_nb_n)^T$ and $a\oslash b\!=\!(a_1/b_1, \cdots, a_n/b_n)^T$.
 The notation $\mathbf 1_n$ denotes an $n$-dimensional vector with each entry being $1$. We use $\mathbb{R}_{+}^n$ ($\mathbb{R}_{++}^n$) to denote the set of $n$-dimensional vectors with nonnegative (positive) components.
 For vectors $\a_1, \cdots, \a_m$, the notation ${\rm col}\{\a_i\}_{i=1}^m$ denotes a new vector $(\a_1^T, \cdots, \a_m^T)^T$. For matrices $M_1, \cdots, M_m$, we use ${\rm diag}\{M_i\}_{i=1}^m$ to denote block diagonal matrix with $i$-th block being $M_i$. The inner product between matrices is denoted as $A:B\!=\!{\rm tr}(A^TB)\!=\!\sum_{i,j}a_{ij}b_{ij}$, where ${\rm tr}$ denotes the matrix trace.
 For a map $g: \mathbb{R}^n\!\ra\! \mathbb{R}^m$ which is differentiable at $x$, we use $\n g(x)$ to denote
the matrix whose rows are the gradients of the corresponding entries in the vector  $g(x)$.

{\bf B. Graph Theory.} Consider a graph
$\mathcal{G}=(\mathcal{V}, \mathcal{E})$, where $\mathcal{V}=\{1,2,\cdots,N\}$ is the set of nodes
representing $N$ agents and $\mathcal{E} \subset \mathcal{V} \times \mathcal{V}$ is the set of
edges of the graph.
The graph considered in this paper is undirected in the sense that the edges $(i,j)$
and $(j,i)$ in $\mathcal{E}$ are considered to be the same.
The set of neighbors of node $i$
is denoted by $\mathscr{N}_i=\{j\in \mathcal{V}: (j,i)\in \mathcal{E}, j\neq i\}$.
We use the symbol $\cup$ to denote the graph union.
We say that  a collection of graphs is jointly connected if the union of its members is a connected graph.
A collection of $\{\mathcal{G}_1, \cdots, \mathcal{G}_S\}$ is jointly connected if and only if the matrix $\mc{L}_1+\cdots+\mc{L}_S$   has a simple zero eigenvalue, where $\mc{L}_1, \cdots, \mc{L}_S$ are respectively the Laplacians of the graphs  $\mathcal{G}_1, \cdots, \mathcal{G}_S$.

\section{Problem Formulation}
Consider an optimization  problem on a graph $\mc{G}=(\mc{V}, \mc{E})$. Each agent $i \in \mc{V}$ has a local cost function $f_i(x)$  and a group of local inequality constraints  $g_{ij}(x)\leq 0$, $j=1, \cdots, r_i$ and equality constraints $h_{ik}(x)=0, k=1, \cdots, s_i$, where $r_i$ and $s_i$ are nonnegative integers. If there is no constraints  for agent $i$, one simply sets corresponding constraint functions to be zero.   The total cost function of the network is given by sum of all local functions, and  the optimization is to minimize the global cost function of the network while satisfying $N$ group of  local constraints, given explicitly as follows,
\begin{eqnarray}\label{optimization}
\mathcal{P}:
\left\{
\begin{array}{llll}
  {\rm minimize}   && \tilde f(x)=\sum_{i=1}^{N}f_i(x), \label{objective}\\
  {\rm subject~to} &&  g_i(x)\preccurlyeq 0,  \label{constraint}\\
  {\rm }           && h_i(x)=0, i=1, \cdots, N. \label{constraint2}
\end{array}
\right.
\end{eqnarray}
where  $f_i: \mathbb{R}^n \rightarrow \mathbb{R}$,
$g_i=(g_{i1}, \cdots, g_{ir_i})^T: \mathbb{R}^n \rightarrow \mathbb{R}^{r_i}$,  and
$h_i=(h_{i1}, \cdots, h_{is_i})^T:\mathbb{R}^n \rightarrow \mathbb{R}^{s_i}$  are respectively the local cost, inequality constraint and equality constraint on node $i$.

Let $x^*$ be an optimal solution, if exists,  to the problem \eqref{optimization}.  If additional assumptions on the constraint functions,  called constrained qualifications,    are satisfied,  then the following classical KKT conditions hold at the minimizer  $x^*$: there exist $\lambda_{ij}^*$ and $\nu_{ij}^*$ such that, for $i=1, \cdots, N$,
\begin{subequations} \label{KKT}
    \begin{empheq}[left={\mathcal{K}\mathcal{K}\mathcal{T}:\empheqlbrace\,}]{align}
      & g_{ij}(x^*) \leq 0, j=1, \cdots, r_i, \label{KKTa}\\
      & h_{ij}(x^*) = 0, j=1, \cdots, s_i,\label{KKTb}\\
      & \lambda_{ij}^*\geq 0, j=1, \cdots, r_i, \label{KKTc}\\
      &  \lambda_{ij}^*g_{ij}(x^*)=0, j=1, \cdots, r_i,\label{KKTd}\\
      &\textstyle \sum \nolimits_{i=1}^{N}\! \nabla f_i(x^*)\!+\!\sum \nolimits_{i=1}^{N}\!\sum\nolimits_{j=1}^{r_i}\lambda_{ij}^* \nabla g_{ij}(x^*)\nonumber\\
      & \hspace{1.3cm}+\textstyle\sum \nolimits_{i=1}^{N}\sum \nolimits_{j=1}^{s_i}\nu_{ij}^*\nabla h_{ij}(x^*)=0. \label{KKTe}
    \end{empheq}
\end{subequations}
A widely used constrained qualification is the Slater's constrained qualification (SQC):  there exists $x\in \mathbb{R}^n$ such that $ g_i(x) \!\prec \!0$ and  $h_i(x)\!=\!0$ for $i=1,\cdots, N$.  In other  words, assuming SQC,
``\scalebox{0.8}[1]{$x^*$ solves ($\mc P$)''} $\Rightarrow$ ``\scalebox{0.8}[1]{$\exists$ a set of $(\lambda_{ij}^*, \nu_{ij}^*)$} \scalebox{0.8}[1]{together with $x^*$ solving} ($\mc K\mc K\mc T$)''. Furthermore, for convex problem, this implication is bidirectional. Refer to   \cite{boyd2004} for details.

While SQC ensures the existence of multipliers satisfying
($\mc K\mc K\mc T$), it does not grantee uniqueness.
Closely tied to this direction is the linear independent constraint
qualification (LICQ), which is  stronger than SCQ.
Using
$\n_{\hspace{-0.1cm}J_i} g_i(x^*)$ to denote the submatrix of
$\n g_i(x^*)$ given by rows with indices in
$J_i(x^*)\!=\! \left\{j | g_{ij}\left(x^{*}\right)\!=\!0\right\}$,
the LICQ is defined  as
\begin{align}\label{LICQ}\renewcommand\arraystretch{0.8}
\operatorname{rank}\left[
\begin{array}{lll}
\nabla h_i(x^*) \\
\nabla_{\hspace{-0.1cm}J_i}g_i(x^*)
\end{array}\right]=s_i+|J_i(x^*)|, \quad i=1, \cdots, N,
\end{align}
here $|\cdot|$ denotes the set cardinality.
By assuming LICQ,  one obtains a result \cite{wachsmuth2013} on the existence and uniqueness of multipliers satisfying \eqref{KKT},
\begin{align}\label{uniqueMultiplier}
\scalebox{0.8}[1]{``$x^*$ solves ($\mc P$)''} \Rightarrow \scalebox{0.8}[1]{``$\exists$ a unique $(\lambda_{ij}^*, \nu_{ij}^*)$ together with $x^*$ solving ($\mc K\mc K\mc T$)''}.
\end{align}
This direction is bidirectional if the problem is convex.  In what follows, we assume that $f_i, g_i, i=1, \cdots, N$  are convex and twice differentiable, and $h_i, i=1, \cdots, N$ are affine, so that  the optimization problem \eqref{optimization} is convex.
Also, we assume the existence of optimal solutions without giving explicit conditions due to space limitation; interested readers can refer to \cite{boyd2004}.
We further assume strict convexity on at least one function among  $\{f_1, \cdots, f_i\}$,  so that there exists
at most one global optimal solution to the problem \eqref{optimization}.
With these assumptions, the unique optimal solution  is denoted by $x^*\in \mathbb{R}^n$  and the corresponding optimal value is denoted by $p^*=\tilde f(x^*)$.

\section{Distributed Optimization Dynamics}

We use $N$  agents $\{1, \cdots, N\}$ to solve the convex optimization problem \eqref{optimization} in a distributed way. Each agent, say $i$,  is a dynamical system with states  $(x_i(t), \theta_i(t), \lambda_i(t), \nu_i(t))$  and communicates its estimate $x_i$ of the optimal solution $x^*$ to its neighboring agents. We consider the general case that the communication is corrupted by noises which lie in a filtered probability space $(\Omega, \mathcal{F}, \mathbb{P}, \{\mathcal{F}_t\}_{t\geq 0})$, where $\{\mathcal{F}_t\}_{t\geq 0}$ is a sequence of increasing $\sigma$-algebras with $\mathcal{F}_{\infty}\subset \mathcal{F}$ and $\mathcal{F}_0$ containing all the $\mathbb{P}$-null sets in $\mathcal{F}$. In more details, for any two neighbor agents $i,j$ with communication channel $(j,i)$ connecting them, the ideal relative information  $(x_j\!-\!x_i)$  transmitted  in this channel is corrupted by state dependent noise $\sigma_{ji}\xi_{ji}(x_j-x_i)$ with  $\xi_{ji}(t)\in \mathbb{R}$ being independent standard white noises adapted to the filtration $\{\mathcal{F}_t|t\geq 0\}$ and $\sigma_{ji}\geq 0$ being noise intensity.
This kind of noise indicates  that the closer the agents are to each other, the smaller the noise intensities.
The noise of this type is multiplicative in nature.  While additive noise provides a natural intuition, multiplicative noise model also has its practical background.
For example, it can model the impact of  quantization error, as well as the effect of a fast-fading communication channel. Also, lossy communication induced noises and imperfect sample induced noises are all multiplicative noises. We also note that the issue of adopting multiplicative noise in inter-agent communication has been well explained in existing literatures (see for example references \cite{ni2013}, \cite[Remarks 1-2]{li2014}, \cite{carli2008}, \cite{wang2013}  and references therein).
By adopting this noise model, we design for each agent $i\in \{1, \cdots, N\}$ the following optimization dynamics
\begin{subequations} \label{algorithm}
 \begin{empheq}[left={\empheqlbrace\,}]{align}
  &\dot x_i=c\sum_{j\in \mathscr{N}_i(t)}(1+\sigma_{ji}\xi_{ji})(x_j-x_i)-\!\nabla f_i(x_i) -\theta_i\nonumber\\
  &\hspace{1.5cm}-\!\sum\nolimits_{j}^{r_i}\lambda_{ij}\nabla g_{ij}(x_i)-\sum\nolimits_{j}^{s_i}\nu_{ij}\nabla h_{ij}(x_i), \label{algorithma}\\
 & \dot \theta_i=-c\sum_{j\in \mathscr{N}_i(t)}(1\!+\!\sigma_{ji}\xi_{ji})(x_j-x_i) \label{add1}\\
  &\dot \lambda_{ij}=\frac{\lambda_{ij}}{1+\eta_{ij}\lambda_{ij}}  g_{ij}(x_i), \hspace{0.5cm}j=1, \cdots, r_i, \label{algorithmb}\\
  &\dot \nu_{ij}=h_{ij}(x_i),\hspace{2cm} j=1, \cdots, s_i, \label{algorithmc}
    \end{empheq}
\end{subequations}
where $x_i, \theta_i \in \mathbb{R}^n$, $\lambda_{ij}, \nu_{ij}\in \mathbb{R}$,  $\eta_{ij}$ are positive parameters and $c>0$ is the coupling strength.

The equations  \eqref{algorithma}-\eqref{add1} are motivated by  \cite[Eq. (3)]{kia2015}, which however  does not consider optimization constraints, communication noises, and more importantly time-varying network. This series of hard problems are tackled in our paper.
The equation \eqref{algorithmb} is motivated by \cite[Eq. (4)]{brunner2012}, but differs from {\rm \cite{brunner2012}} since it is distributed by using $N$ agents on a network and considers the communication noises among agents and the stochastic networks.
The dynamics \eqref{algorithmc} can be obtained  by maximizing  $\Phi$ in \eqref{barlagrange}  with respect to $\nu_{ij}$  via gradient ascent $\dot {\nu}_{ij}=\n_{\nu_{ij}} \Phi$.

The algorithm \eqref{algorithm} is designed for time-varying  networks.
Assume  that there are $S\! \in\! \mathbb N$ possible graphs $\{\mathcal{G}_1, \cdots, \mathcal{G}_S\}$, among which the network structure is switched. We use a continuous-time Markov chain $\sigma: [0, \infty)\!\rightarrow\! \mathbb S\!:=\!\{1, \cdots, S\}$  to describe this switching.
To analyze the stability of equation \eqref{algorithm}, we rewrite it into a switching SDE as (refer to \hyperlink{appendix1}{Appendix A}  for detailed derivation and implicit definition of symbols in the  equation below),
\begin{subequations} \label{X}
 \begin{empheq}[left={\hspace{-0.2cm}\empheqlbrace\,}]{align}
 &d\hx\!=\![-c\boldsymbol{\mathcal{L}}_{\sigma(t)}\hx\!-\!\htheta\!-\! \nabla F(\hx)\!-\!\boldsymbol{\lambda} \!\odot\! \nabla G(\hx) \!-\!\boldsymbol{\nu} \!\odot\! \nabla H(\hx)]dt\!+\!c\mathcal{M}_{\sigma(t)}d\hw, \label{Xa}\\
  & d\htheta=c\boldsymbol{\mathcal{L}}_{\sigma(t)}\hx dt-c\mathcal{M}_{\sigma(t)}d\hw, \label{Xb} \\
 &d \boldsymbol{\lambda} =\left[\boldsymbol{\lambda} \oslash (\mathbf{1}+\boldsymbol{\eta} \odot \boldsymbol{\lambda})\right] \odot  G(\hx)dt, \label{Xc} \\
 &d \boldsymbol{\nu}= H(\hx)dt. \label{Xd}
 \end{empheq}
\end{subequations}
Note that we have assumed the existence and uniqueness of an optimal solution $x^*$. Defining $\hx^*=\mathbf 1_N \otimes x^*$,   it follows from \eqref{uniqueMultiplier} that there is a unique pair of $(\hlambda^*, \hnu^*)$ satisfying \eqref{KKT}, where $\hlambda^*$ and $\hnu^*$ are respectively stacked vectors of $\lambda_{ij}^*$ and $\nu_{ij}^*$ in \eqref{uniqueMultiplier}.
These give rise to another vector $\htheta^*$ defined by
\begin{align}\label{conditionTheta}
\hspace{-0.35cm}\htheta^*\!+\! \nabla F(\hx^*) \!+\!\boldsymbol{\lambda}^* \!\odot\! \nabla G(\hx^*)\!+\!\boldsymbol{\nu}^*\! \odot\! \nabla H(\hx^*)=0,
\end{align}
Obviously,  $(\hx^*, \htheta^*, \hlambda^*, \hnu^*)$ satisfies \eqref{X} and consequently  is an equilibrium of \eqref{X} (note that stochastic disturbances vanish at this equilibrium).


In the rest of this paper, we will utilize the averaging method used in \cite{ni2016a,ni2016b} to analyzed the stability of the equilibrium  $(\hx^*, \htheta^*, \hlambda^*, \hnu^*)$ for  \eqref{X}. The following assumptions are made:
\begin{assumption}\label{A1}
The LICQ  defined in \eqref{LICQ} is satisfied.
\end{assumption}
\begin{assumption}\label{A2}
The noise intensity $\sigma_{ij}$ has an upper bound
$\sigma_{ij}\leq \kappa$ for some positive constant $\kappa$, where $i, j=1, \cdots, N$.
\end{assumption}

\section{Distributed Optimization Under Fixed Network}
In this section, we assume that the network is  fixed, so that the time-varying  graph Laplacian $\boldsymbol{\mathcal{L}}_{\sigma(t)}$ and  the diffusion term $\mathcal{M}_{\sigma(t)}$ in the dynanmics \eqref{X}  are replaced  with fixed ones $\boldsymbol{\mathcal{L}}$ and $\mathcal{M}$,  respectively.  The following assumption is made:
\begin{assumption}\label{A3}
The coupling strength $c$ for fixed network satisfies $0<c<\frac{2}{3}\kappa^{-2}$.
\end{assumption}

The following theorem states that the trajectory of our optimization dynamics converges to its equilibrium $(\hx^*, \htheta^*, \hlambda^*, \hnu^*)$ under appropriate conditions. While this result is obtained for fixed network, the multi-scale method developed in next section can transform the analysis of the optimization algorithm under switching networks to that under fixed one.
\begin{theorem}\label{theorem3}
For the constrained optimization problem \eqref{optimization}, suppose there are $N$ agents  with each agent dynamics  given by \eqref{algorithm},  and they are connected across a fixed network (the graph Laplacian and diffusion term in \eqref{X} are now denoted as $\boldsymbol{\mathcal{L}}$ and $\mathcal{M}$ respectively). Under Assumptions \ref{A1}-\ref{A3} with $\kappa \leq \sqrt{\lambda^{\mc L}_2}/2$ with  $\lambda^{\mc L}_2$ the smallest nonzero eigenvalue of the graph Laplacian $\mc L$,   then for any trajectory of \eqref{X} with initial conditions
$\lambda_{ij}(0)> 0$ and $\sum_{i=1}^N\theta_i(0)=0$, one has  $\lim_{t\rightarrow \infty}\|x_i(t)-x^*\|=0$ for $i=1, \cdots, N$ almost surely.
\end{theorem}

{\it Proof}  We use the method of Lyapunov function to prove the stability, with the key to construct an appropriate Lyapunov function. Here we only give a proof skeleton, with more details being put in Appendices B, C, and D.

${\rm \mathbf 1^{o}}$ Construct a Lyapunov candidate  $V(\hx, \htheta, \hlambda, \hnu)=V_1+V_2+V_3+V_4$ with
\begin{align*}
&  V_1=\frac12 \|\hx-\hx^*\|^2+\frac 1 2 \|(\hx-\hx^*)+(\htheta-\htheta^*)\|^2,\\
&  V_2=\frac12 \sum\nolimits_{i=1}^{N} \sum\nolimits_{j=1}^{r_i}\eta_{ij}(\lambda_{ij}-\lambda_{ij}^*)^2,\\
&  V_3=\sum\nolimits_{(i, j)\in \Omega}D_{\phi}(\lambda_{ij}, \lambda_{ij}^*) +\sum\nolimits_{(i, j)\notin \Omega}(\lambda_{ij}-\lambda_{ij}^*)^2,\\
 & V_4=\frac12 \sum\nolimits_{i=1}^{N} \sum\nolimits_{j=1}^{s_i}(\nu_{ij}-\nu_{ij}^*)^2,
\end{align*}
where $\Omega\!=\!\{(i,j)| \lambda_{ij}^* \!\neq\! 0\}$ and $D_{\phi}(\lambda_{ij}, \lambda_{ij}^*)\geq 0$ is the Bregman divergence between $\lambda_{ij}$ and $\lambda_{ij}^*$ with respect to $\phi(x)=x\ln x$ (refer to \cite{bregman1967} for the definition of Bregman divergence).   Therefore, $V\!\geq \!0$,  and also $V(\hx, \htheta,  \hlambda,  \hnu)=0$ if and only if $(\hx, \htheta, \hlambda, \hnu)\!=\!(\hx^*,  \htheta^*, \hlambda^*,  \hnu^*)$.

${\rm \mathbf 2^{o}}$ We now calculate the action of the operator $\mathcal{A}$  on  the function $V$ (Recall that, for an SDE $dx=a dt+bdw$,  $\mc{A}V=\n V \cdot a+\frac{1}{2}{\rm tr}(b^T\n^2Vb)$).
Defining an Lagrangian $ \Phi: \mathbb{R}^{nN} \times \mathbb{R}^{nN}\times \mathbb{R}^r_+ \times \mathbb{R}^s\rightarrow \mathbb{R}$ as
\begin{align}\label{barlagrange}
\hspace{-0.3cm} \Phi(\hx,  \htheta, \hlambda,  \hnu)\!=\! \Psi \!+\!\sum_{i=1}^{N}\!\sum_{j=1}^{r_i}\!\lambda_{ij}g_{ij}(x_i)\!+\!
\!\sum_{i=1}^{N}\!\sum_{j=1}^{s_i}\!\nu_{ij}h_{ij}(x_i)
\end{align}
with
 \begin{align}\label{psi}
\Psi(\hx, \htheta)=\sum_{i=1}^{N}f_i(x_i)+(\hx-\hx^*)^T\htheta+\hbar\hx^T\hL \hx, \quad  \hbar=\frac{1}{2}(c-\frac{1}{2}c^2\kappa^2),
 \end{align}
 we can show in \hyperlink{appendix3}{Appendix B} that
\begin{align}\label{e:appendix}
\mathcal{A}V
& \leq  \Phi(\hx^*\hspace{-0.12cm}, \htheta, \hlambda, \hnu)- \Phi(\hx, \htheta^*\hspace{-0.15cm}, \hlambda^*, \hnu^*)\!-\!(\hbar \lambda^{\mathcal G}_2\!-\!1)(\hx\!-\!\hx^*)^T(\hx\!-\!\hx^*),
\end{align}
where $\lambda^{\mathcal G}_2$ is the smallest nonzero eigenvalue of the connected graph $\mathcal G$.
Noting that $(x^*, \hlambda^*, \hnu^*)$
satisfies \eqref{KKT},
we show in \hyperlink{appendix4}{Appendix C} that the following saddle point condition holds
\begin{align}\label{SaddlePointCondition}
 \Phi(\hx^*, \htheta, \hlambda, \hnu)
\!\leq\!
 \Phi(\hx^*, \htheta^*, \hlambda^*, \hnu^*)
\leq
 \Phi(\hx, \htheta^*, \hlambda^*, \hnu^*).
\end{align}
 Therefore,   $\mathcal{A}V\!\leq\! 0$ if $\hbar \lambda^{\mathcal G}_2\!-\!1>0$, which is guaranteed by $\kappa \leq \scalebox{0.8}[0.7]{$\sqrt{\lambda_2^{\mc G}}$}/2$.

${\rm \mathbf 3^{o}}$
Letting $\mathcal{A}V=0$, we show in \hyperlink{appendix4}{Appendix D}  that $(\hx, \htheta, \hlambda, \hnu)=(\hx^*, \htheta^*, \hlambda^*, \hnu^*)$. The application of the stochastic version of the Lasalle's invariance principle in \cite{mao1999} yields that the equilibrium $(\hx^*, \htheta^*, \hlambda^*, \hnu^*)$  of system \eqref{algorithm}  is asymptotically stable almost surely.
\qed

\section{Multi-Scale Analysis of Distributed Optimization Dynamics}
This section analyze the optimization dynamics  under  switching networks by proposing a multi-scale analysis. This method can reduce the analysis of optimization algorithm under switching topologies to  that under a fixed one.

The switching feature encoded in $\sigma(t)$ makes difficult the convergence analysis of system \eqref{X}.
To cope with this difficulty, we adopt an idea of  SA proposed in our earlier works \cite{ni2016a,ni2016b}.
The basic idea is to approximate in an appropriate  sense  the switching system \eqref{X} using a non-switching system,  called the average system.   A detailed construction of such an average system is provided in this section by resorting to the multi-scale analysis.
Noting that adding or deleting even one edge in the graph may result in the change of the network structure, and also noting that the graph of the network is large scale (i.e. the number of nodes and the number of edges are extremely large), the change of the network structure takes place more readily than the evolution of  the optimization dynamics on the nodes. That is, the state $\s(t)$ of the network structure changes faster than the state ${\rm col}[\hx(t), \htheta(t), \hlambda(t), \hnu(t)]$ of the optimization dynamics. We use a small parameter $\a>0$ to resale $\s(t)$ so that $\s(t/\a)$ is a fast process and ${\rm col}[\hx(t), \htheta(t), \hlambda(t), \hnu(t)]$ is a slow process.
Correspondingly, the stochastic differential equation \eqref{X} under the above re-scaling
admits the following form
\begin{subequations} \label{Xrescale}
 \begin{empheq}[left={\empheqlbrace\,}]{align}
 &d\hx^{\alpha}\!=\![-c\boldsymbol{\mathcal{L}}_{\sigma(\hspace{-0.05cm}\frac{t}{\alpha}\hspace{-0.05cm})}\hx^{\alpha}-\htheta^{\alpha}
 -\nabla F(\hx^{\alpha}\!)-\hlambda^{\alpha} \odot \nabla G(\hx^{\alpha})\\
& \hspace{2.2cm}\!-\hnu^{\alpha} \odot \nabla H(\hx^{\alpha})]dt+c\mathcal{M}_{\sigma(\frac{t}{\alpha})}d\hw(t),\\
&  d\htheta^{\alpha}\!=c\boldsymbol{\mathcal{L}}_{\sigma(\frac{t}{\alpha})}\hx^{\alpha},\\
&  d \hlambda^{\alpha}=\left[\hlambda^{\alpha} \oslash (\mathbf{1}+\eta \odot \hlambda^{\alpha})\right] \odot  G(\hx^{\alpha}),\\
& d \hnu^{\alpha}=  H(\hx^{\alpha}).
 \end{empheq}
\end{subequations}
The time re-scaling  is crucial in the method of SA and its role  will be seen later.
Some assumptions on the Markov chain $\sigma(t)$ are now given below.


\subsection{Assumptions on The Markov Switching Network.}
It is well known that the statistics of a Markov chain defined over
$\mathbb{S}=\{1, \cdots, S\}$  is identified by
an initial probability distribution $\pi_0=[\pi_{01}, \cdots, \pi_{0S}]^T$  with $\pi_{0i}=\mathbb{P}(\sigma(0)=i)$  and by a Metzler
matrix \scalebox{0.9}[1]{$\mathcal {Q}=(q_{ij})_{S \times S} \in \mathbb{R}^{S \times S}$}. This matrix is also called the infinitesimal
generator of the Markov chain and it describes the transition probability as
\begin{align*}
  \mathbb{P}\{\sigma(t+\Delta t)=j|\sigma(t)=i\}
  =  \left\{
  \begin{array}{llllll}
    q_{ij}\Delta t+o(\Delta t), &   i \neq j,\\
   1+q_{ii} \Delta t+o(\Delta t), &   i = j,
  \end{array}
  \right.
\end{align*}
where $\Delta t>0$,  $q_{ij}\geq 0$ ($i \neq j$)
 is the transition rate from state $i\in \mathbb{S}$ at time $t$ to state
 $j\in \mathbb{S}$ at time $t+\Delta t$,
 and $q_{ii}=-\sum_{j\neq i}q_{ij}$, $\lim_{\Delta t\rightarrow \infty}o(\Delta t)/\Delta t=0$ (Note that $|q_{ij}|<\infty$ since $\sigma(t)$ is a finite state Markov chain \cite[pp.150-151]{freedman1983}).
More specifically,  at time $t$ the state of the Markov chain is determined according to the probability
distribution $\pi(t)=(\pi_1(t), \cdots, \pi_S(t))$ with $\pi_s(t)$ being the probability that at time $t$ the Markov system is in the state $s \in \mathbb{S}$. The normalization
condition $\sum_{s=1}^{S}\pi_s(t)=1$ is usually assumed.
Letting $\Delta t\rightarrow 0$, the infinitesimal form of the Markov dynamics reads as  $\dot\pi_s(t)=\sum_{i=}^{S}\pi_i(t) q_{is}$, $s=1, \cdots, S$. In a compact form,  the distribution $\pi(t)$ for $\sigma(t)$ obeys the differential equation
\begin{align}\label{dynamicsPit}
\dot \pi(t)=\mathcal {Q}^T \pi(t).
\end{align}
Since we are interested in the asymptotic behavior of the system, we will assume  that the  probability distribution $\pi(t)$ of $\sigma(t)$ is stationary, and it is denoted by $\pi$, which  is defined as
\begin{align}\label{StationaryDistribution}
 \mathcal {Q}^T  \pi=0, \quad \sum \nolimits_{i=1}^{S} \pi_i=1, \quad \pi_i>0.
\end{align}
The existence of the stationary distribution $\pi$ satisfying \eqref{StationaryDistribution} can be guaranteed by the ergodicity of $\sigma$, namely, all graphs $\{\mathcal{G}_1, \cdots, \mathcal{G}_S\}$ can be visited infinitely  often under the switching $\s$.
The joint connectivity of the network is also needed for later analysis and it is also assumed here.
\begin{assumption}\label{AssumptionSwitching}
The finite-state Markov process $\sigma(t)$ describing the switching networks has a stationary probability
distribution $\pi=(\pi_1, \cdots, \pi_N)^T$ satisfying (\ref{StationaryDistribution}), and the
union graph $\cup_{s\in \mathbb{S}}\mathcal{G}_s$ is connected.
\end{assumption}

Due to Assumption \ref{AssumptionSwitching}, the eigenvalues $\bar \lambda_1\!<\!\bar\lambda_2\!\leq\! \cdots\leq \bar \lambda_N$ of the matrix $\mc{\bar L}:=\sum_{i=1}^S\!\mathcal{L}_i$  has a simple zero eigenvalue $\bar \lambda_1\!=\!0$.
Denote $\pi_{\text{min}}=\min\{\pi_1, \cdots, \pi_S\}$ and $\pi_{\text{max}}=\max\{\pi_1, \cdots, \pi_S\}$. For later use, we modify Assumption \ref{A3} in fixed topology to generalize to switching topology as follows.

\begin{assbis}\label{AA3}
The coupling strength $c$ under the switching network satisfies $0<c< \scalebox{0.95}[0.92]{$\frac{2}{3}(\pi_{\text{min}}/\pi_{\text{max}})$}\kappa^{-2}$.
\end{assbis}

The following  lemma is a slight modification of \cite[Lemma 4.2]{fragoso2005}, and it can be proved similarly as in \cite{fragoso2005}.
\begin{lemma} \label{LemmaDifferentialMarkov}
  Suppose that  $V(t)$ is $\mathcal{F}$-measurable and  $\mathbb{E}\{V(t)\mathbf{1}_{\{\sigma(t)\!=\!i\}}\}$ exists,  where $\mathbf{1}_{\{\sigma(t)=i\}}$ is the indicator function of the event $\{\sigma(t)=i\}$. Then  $\mathbb{E}[V(t) \cdot d(\mathbf{1}_{\{\sigma(t/\a)=s\}})]=\frac{1}{\alpha}\sum_{j=1}^S q_{js}\mb{E}[V(t)]dt+o(dt)$ holds for each $s \in  \mathbb S$.
\end{lemma}

\subsection{Stochastic Averaging Method for Switching Networks.}
The properties of the solutions to the SDE \eqref{Xrescale} is included in its backward Kolmogorov equation which is a partial differential equation determined by  the infinitesimal of \eqref{Xrescale}. Although the analytic solutions to the backward Kolmogorov equation are hard to obtain, we focus on those solutions which are Taylor series in term of $\alpha$. We use the first term in the series as an approximate solution. It will be shown that this approximate solution satisfies another backward Kolmogorov equation, which is called as the averaged backward Kolmogorov equation. This averaged backward Kolmogorov equation is nothing but the one whose operator is the  average of infinitesimals of all the subsystems in the switched system \eqref{Xrescale}. Corresponding to this average backward Kolmogorov equation, there is an SDE which is time-invariant and is called the average SDE. The analysis of the original SDE \eqref{Xrescale} can approximately be transformed to  the average SDE. Therefore, the problem in the switching case in this section can be reduced to the one in fixed case in last section.

{\bf Backward Kolmogorov equation for SDE \eqref{Xrescale}.}
Denote $Z$ as the state of \eqref{Xrescale}; that is, $Z\!=\!{\rm col}[\hx^{\alpha}, \htheta^{\alpha}, \hlambda^{\alpha}, \hnu^{\alpha}]$.
The stochastic process $(Z(t), \sigma(t))$, rather than the $Z(t)$,  is a Markovian process, whose
infinitesimal is $\mathcal{A}^{\alpha}=\scalebox{1}[0.96]{$\frac{1}{\alpha}$}\mathcal{Q}+\mathcal{A}_s$, where $\mathcal{A}_s$ is the
infinitesimal of the $s$-subsystem of \eqref{X}, $\mathcal{Q}$ is the
infinitesimal of the Markov chain $\sigma(t)$.
This means that the average number of jumps of $\s(t)$ per unit of time is proportional to $1/\a$.
 Let $\phi$ be a sufficiently smooth real-valued function defined on the state space $(Z, \sigma)$
 and let \scalebox{0.9}[1]{$W(t, Z, s)=\mathbb{E}\left[\phi(Z(t), \sigma(t))|Z(0)=Z,\sigma(0)=s\right]$}.
From the standard analysis in stochastic theory, \scalebox{0.9}[1]{$W(t, Z, s)$} is a unique bounded classical solution to the following partial differential equation with the initial data $W(0, Z, s)=\phi(Z, s)$:
\begin{align}\label{BackKol}
  \frac{\partial}{\partial t}W(t, Z, s)=\frac{1}{\alpha} \mathcal{Q}  \overrightharpoon{W}(t, Z)[s]+\mathcal{A}_sW(t, Z, s),
\end{align}
where $\overrightharpoon{W}(t, Z)\!=\!\big(W(t, Z, 1), \cdots, W(t, Z, S) \big)^T$ and $\mathcal{Q}\overrightharpoon{W}(t, Z)[s]$
denotes the $s$-row of the matrix $\mathcal{Q}\overrightharpoon{W}(t, Z)$.
The partial differential equation \eqref{BackKol} is termed  as  the backward Kolmogorov equation associated with
the SDE \eqref{Xrescale}.

\scalebox{0.9}[1]{\bf Average backward Kolmogorov equation for SDE \eqref{Xrescale}.}
We single out for \eqref{BackKol} an approximate solution of the form  $W=W_0+\a W_1+\mathcal {O}(\a^2)$.
Inserting this expression  into \eqref{BackKol} and equating coefficients of $\alpha^{-1}$ on both sides yields $\mathcal{Q}\overrightharpoon{W}_0=0$.
Due to  Assumption \ref{AssumptionSwitching} or its equivalent characterization of \eqref{StationaryDistribution},  one sees that the null space of the adjoint generator $\mathcal{Q}^T$ consists of only constants, which also amounts to saying that the null space of its  infinitesimal generator $\mathcal{Q}$ consists of only constant functions. This fact, together with $\mathcal{Q}\overrightharpoon{W}_0=0$, implies that $\overrightharpoon W_0$ is a function independent of the switching mode $s\in \mathbb{S}$; that is, $W_0(t,Z,1)=W_0(t,Z,2)=\cdots=W_0(t,Z,S)$.     For ease of notation, we denote them by $W_0(t, Z)$.
Similarly, inserting the expression of $W$  into \eqref{BackKol} and equating coefficients of $\alpha^{0}$ on both sides yields
\begin{align}\label{QW}
  \mathcal{Q}\overrightharpoon W_1=
  \left(
    \begin{array}{c}
      \frac{\partial W_0}{\partial t}-\mathcal{A}_1W_0 \\
      \vdots \\
      \frac{\partial W_0}{\partial t}-\mathcal{A}_SW_0 \\
    \end{array}
  \right).
\end{align}
Recall the Freddhom alterative (see e.g. \cite[pp.641 Theorem 5(iii)]{evans2010} for general case, but only a simple version in finite dimension is needed here) which deals with  the solvability of an  inhomogeneous linear algebraic equation $Ax=b$ with $A\in \mb{R}^{n \times n}$ and $b\in \mb{R}^n$; that is,  this inhomogeneous equation is solvable if and only if  $b$ belongs to the column space of $A$, which is the orthogonal complement of ${\rm ker}(A^T)$.
A direct application of  the Freddhom alterative to the  equation \eqref{QW} tells us  that
$(\frac{\partial W_0}{\partial t}-\mathcal{A}_1W_0, \cdots, \frac{\partial W_0}{\partial t}-\mathcal{A}_SW_0)^T$  is perpendicular to the null space of $\mathcal{Q}^T$. Noting that ${\rm Null}(\mathcal{Q}^T)=\pi$ in \eqref{StationaryDistribution},
it is nature to have $\sum_{s=1}^{S}\pi_s \left(\frac{\partial W_0}{\partial t}-\mathcal{A}_sW_0 \right)=0$, which gives rise to an \emph{average backward equation}
\begin{align}\label{averageBK}
\frac{\partial W_0}{\partial t}=\mathcal{A}_{\pi}\, W_0,
\end{align}
where $\mathcal{A}_{\pi}=\sum_{s=1}^{S}\pi_s \mathcal{A}_s$ is the stochastic average of the infinitesimal generators $\mathcal{A}_s$ with respect to the invariant measure $\pi$.
As an conclusion, we have  shown that  the first term $W_0$ in the series of $W$ is a solution to another backward equation specified by the average operator $\mathcal{A}_{\pi}$.

{\bf Average SDE for \eqref{Xrescale}.}
We proceed to construct for the average backward equation \eqref{averageBK} an SDE whose  infinitesimal is exactly $\mathcal{A}_{\pi}$. We call such an SDE as the average equation for \eqref{Xrescale}.
Denoting the drift vector  for the $s$-subsystem in equation  \eqref{Xrescale} by  $\boldsymbol{a}_s$  and the diffusion matrix (the diffusion matrix for an SDE $dx\!=\!fdt\!+\!gdw$ is defined as $gg^T$) by $\boldsymbol{\Gamma}_s$, and noting
$\mathcal{A}_s\!=\!\boldsymbol{a}_s \!\cdot\! \boldsymbol{\nabla}+\frac{1}{2} \boldsymbol{\Gamma}_s\!:\!\boldsymbol{\nabla}\boldsymbol{\nabla}$,
the average operator $\mathcal{A}_{\pi}$ can be calculated as $\mathcal{A}_{\pi}\!=\!\boldsymbol{a}_{\pi} \!\cdot\! \boldsymbol{\nabla}+\frac{1}{2} \boldsymbol{\Gamma}_{\pi}\!:\!\boldsymbol{\nabla}\boldsymbol{\nabla}$,
where
${\boldsymbol{a}}_{\pi}\!=\!\sum_{s=1}^{S}\pi_s \boldsymbol{a}_s$ and
$\boldsymbol{\Gamma}_{\pi}\!=\!\sum_{s=1}^{S}\pi_s \boldsymbol{\Gamma_s}$ are respectively the averages of the drift vector and  the diffusion matrix.
Corresponding to $\mathcal{A}_{\pi}$  above, one can construct an average SDE for \eqref{Xrescale} as follows
\begin{align}\label{Xaverage}
\begin{split}
\left\{
\begin{array}{llllll}
 d\bar \hx\!=\!\big[-\!c\boldsymbol{\mathcal{L}}_{\pi} \bar  \hx\! -\!\bar  \htheta\!-\! \nabla F(\bar \hx)\!-\!\bar\hlambda \!\odot\! \nabla G(\bar \hx)\!-\!\bar\hnu \!\odot \!\nabla H(\bar \hx)\big]dt+c\bar{\mathcal{M}}d\hw(t),\\
 d \bar \htheta=c\boldsymbol{\mathcal{L}}_{\pi} \bar \hx,\\
 d \bar\hlambda=\left[\bar\hlambda \oslash (\mathbf{1}+\eta \odot \bar\hlambda)\right] \odot G(\bar \hx),\\
 d \bar\hnu=  H(\bar \hx),
\end{array}
\right.
 \end{split}
\end{align}
where  $\boldsymbol{\mathcal{L}}_{\pi}=\sum_{s=1}^{S}\pi_s \boldsymbol{\mathcal{L}}_s$ is the stochastic average of the Laplacians $\{\boldsymbol{\mathcal{L}}_s, s\in \mb{S}\}$ for the graphs $\{\mc{G}_s, s\in \mb{S}\}$,  $\bar{\mathcal{M}}\in \mb{R}^{(nN)\times N^2}$ is chosen such that $\bar{\mathcal{M}}\bar{\mathcal{M}}^T=\sum_{s=1}^{S}\pi_s \mc{M}_s\mc{M}_s^T \in \mb{R}^{(nN)\times (nN)}$.

Due to Assumption \ref{AssumptionSwitching} and  similar as the proof of Lemma 3.4 in \cite{ni2012}, the average graph Laplacian $\boldsymbol{\mathcal{L}}_{\pi}$ can be shown to have a simple zero eigenvalue, and thus it can be view as a  Laplacian for a certain fixed connected graph.
Replacing $\boldsymbol{\mc{L}}$ and $\mc{M}$ in Theorem \ref{theorem3} with $\boldsymbol{\mc{L}}_{\pi}$ and $\bar{\mc{M}}$ respectively, modifying Assumption \ref{A3} into Assumption \ref{AA3},  and arguing in a similar line of the proof for Theorem \ref{theorem3}, we  establish a stability result for the average system \eqref{Xaverage}, which will be applied to the stability analysis for the original system \eqref{Xrescale}.

\begin{theorem}\label{theorem4}
Consider the average system \eqref{Xaverage} for the constrained optimization problem \eqref{optimization} under Assumptions \ref{A1}, \ref{A2}, \ref{AA3} and \ref{AssumptionSwitching} with $\kappa \leq \sqrt{\bar \lambda_2}/2$.  Then for any initial conditions with
$\lambda_{ij}(0)> 0$ and $\sum_{i=1}^N\theta_i(0)=0$, one has $\lim_{t\rightarrow \infty}\|\bar x_i(t)- x^*\|=0$ almost surely, where $j=1,\cdots, r_i, i=1, \cdots, N $.
\end{theorem}

The relationship between the solutions of the average system \eqref{Xaverage} and the original system \eqref{Xrescale} is now clarified.
Associate with the systems \eqref{Xrescale} and \eqref{Xaverage}, there are respectively a backward Kolmogorov equation \eqref{BackKol} and an average backward equation \eqref{averageBK}, whose solutions are respectively $W(t, Z, s)$ and $W_0(t, \bar Z)$ (Similar to $W(t, Z, s)$,   $W_0(t, \bar Z)=\mathbb{E}\left[\phi((\bar Z(t), s)|\bar Z(0)=\bar Z)\right]$ by definition).
Since $W=W_0+\alpha W_1+ \mathcal{O}(\alpha^2)$, the solution $W$ has a limit $W_0$ as $\alpha \rightarrow 0$; that is,
$\mathbb{E}\left[\phi(Z^{\a}(t), \sigma(t))|Z(0)=X,\sigma(0)=s\right] \stackrel{\alpha \rightarrow 0}{\longrightarrow} \mathbb{E}\left[\phi((\bar Z(t), s)|\bar Z(0)=\bar Z)\right], \forall t$,
then by \cite[Lemma 4]{ni2016b},  $Z^{\a}(t)$ converges weakly to $ \bar Z(t)$ as $\alpha \rightarrow 0$.
Since $\bar Z(t)$ converges asymptotically to the optimal solution $(\hx^*,(\htheta^*, \hlambda^*, \hnu^*))$ almost surely, it is thus that $ Z^{\a}(t)$ also converges asymptotically to the optimal solution $(\hx^*,(\htheta^*, \hlambda^*, \hnu^*))$.

\begin{theorem}\label{theorem5}
Consider the distributed optimization algorithm \eqref{algorithm}
for the constrained optimization problem \eqref{optimization} under  Assumptions \ref{A1}, \ref{A2}, \ref{AA3} and \ref{AssumptionSwitching} with $\kappa \leq \sqrt{\bar \lambda_2}/2$. Then, with a sufficiently small $\alpha>0$ and    for any initial conditions with
$\lambda_{ij}(0)> 0$ and $\sum_{i=1}^{N}\theta_i(0)=0$, one has $\lim_{t\rightarrow \infty}\|x_i(t)-x^*\|=0$ almost surely, where $j=1,\cdots, r_i, i=1, \cdots, N $.
\end{theorem}

{\it Proof}
The proof  amounts to showing the almost sure stability of the system \eqref{Xrescale}.
To this end, also consider  the Lyapunov candidate $V$ as in Theorem \ref{theorem3} (c.f. \hyperlink{appendix2}{Appendix A2}).
The value of the function $V$ on the trajectory of \eqref{Xrescale} is denoted as $V(t)$.
Also define $V^{[s]}(t)=V(t) \mathbf{1}_{\{\sigma(t/\a)=s\}}$. Obviously, $V(t)=\sum_{s=1}^S V^{[s]}(t)$ almost surely.
Denote $d$ the differential of $V$ along the trajectory of \eqref{Xrescale}.
Therefore, 
\begin{align*}
   \mb{E}[dV]
  &=\Sigma_{s=1}^S \mb{E}[dV^{[s]}]=\Sigma_{s=1}^S \mb{E}[d(V \cdot \mathbf{1}_{\{\sigma(t/\a)=s\}})]\\
  &=\Sigma_{s=1}^S\mathbb{E}[d {V} \cdot \mathbf{1}_{\{\sigma(t/\a)=s\}}]+\Sigma_{s=1}^S\mathbb{E}[V \cdot d(\mathbf{1}_{\{\sigma(t/\a)=s\}})]\\
  &=\Sigma_{s=1}^S \mc{A}_s V \!\cdot\! \mathbb{E}\![\mathbf{1}_{\{\sigma(t/\a)=s\}}]dt\!+\!\Sigma_{s=1}^S\mathbb{E}[V \!\cdot\! d(\mathbf{1}_{\{\sigma(t/\a)=s\}})]\\
  &=\Sigma_{s=1}^S \pi_s\mc{A}_s V dt+\Sigma_{s=1}^S\mathbb{E}[V \cdot d(\mathbf{1}_{\{\sigma(t/\a)=s\}})].
  \end{align*}
By Lemma \ref{LemmaDifferentialMarkov}, the
above \scalebox{0.95}[1]{$\mathbb{E}[d {V} \!\cdot\! \mathbf{1}_{\{\sigma(t/\alpha)=s\}}]$}
can be calculated as
$\mathbb{E}[V(t) \!\cdot\! d(\mathbf{1}_{\{\sigma(t/\a)=s\}})]\!
=\!\frac{1}{\alpha}\Sigma_{j=1}^S q_{js}\mb{E}[V^{[j]}(t)]dt+o(dt).$
Therefore,
\begin{align*}
  \mb{E}[dV]
  &\!=\!\mc{A}_{\pi}V  dt\!+\!(1/\a)\scalebox{0.97}[0.97]{$\sum_{s=1}^S \sum_{j=1}^S$ }q_{js}\mb{E}[V^{[j]}(t)]dt+o(dt)\\
  &\!=\!\mc{A}_{\pi}V  dt\!+\!(1/\a)\scalebox{0.97}[0.97]{$\sum_{s=1}^S \sum_{j=1}^S$ } q_{js}\mb{E}[V^{[j]}(t)]dt+o(dt)\\
  &\!=\!\mc{A}_{\pi}V  dt+o(dt),
  \end{align*}
where the last equality uses $\sum_{s=1}^S q_{js}=0$ for $j\in \{1, \cdots, S\}$.

We now calculate $\mc{A}_{\pi}V$ by arguing in an entirely  similar manner as in the proof of Theorem \ref{theorem3},  with only a modification of  the calculation of $\mc{A}V_1$ in \eqref{AV1} by replacing $\boldsymbol{\mathcal{L}}$ in  \eqref{AV1}-\eqref{MM} with $\boldsymbol{\mathcal{L}}_{\pi}$.  The resulting result on $\mc{A}_{\pi}V$, similar as the one in \eqref{e:appendix},  can be calculated   as
\begin{align*}
\hspace{-0.3cm} \mathcal{A}_{\pi}V
& \!\!\leq\!\! \scalebox{0.95}[1]{$ \Phi_{\pi}(\hx^*\hspace{-0.12cm}, \htheta, \hlambda, \hnu)$}\!-\!\scalebox{0.95}[1]{$ \Phi_{\pi}(\hx, \htheta^*\hspace{-0.15cm}, \hlambda^*, \hnu^*)$}\!-\!(\hbar_{\pi}\bar\lambda_2-1)(\hx\!-\!\hx^*)^T(\hx\!-\!\hx^*),
\end{align*}
where $\Phi_{\pi}(\hx^*\hspace{-0.12cm}, \htheta, \hlambda, \hnu)$ is similarly  defined as $\Phi(\hx^*\hspace{-0.12cm}, \htheta, \hlambda, \hnu)$  by replacing $\hL$  with $\hL_1\!+\!\cdots\!+\!\hL_S$ and $\hbar$ (cf.  \eqref{psi})  with $\hbar_{\pi}=$\scalebox{0.92}[0.95]{$\frac{1}{2}$}$(c\pi_{\text{min}}-\frac{3}{2}c^2\kappa^2\pi_{\text{max}})$.
The rest part of proving asymptotic stability almost surely is similar to the third part of proof for  Theorem \ref{theorem3}. \qed


\section{Simulation}

 Consider the optimization problem \eqref{optimization} on a network with $5$ agents. The five local cost functions for five agents are given as
$f_1(x_1,x_2)=4x_1^2+2x_2$, $f_2(x_1,x_2)=2x_2^2$, $f_3(x_1,x_2)=4x_1$, $f_4(x_1,x_2)=2x_2$, $f_5(x_1,x_2)=e^{3x_1+x_2}$. We assume that agent $1$ has both inequality and equality constraints with constraint functions  $g_1(x_1, x_2)=(x_1-2)^2-x_2+1, h_1(x_1, x_2)=2x_1-x_2$ and agent $2$ has  only inequality constraint with  constraint function   $g_2(x_1,x_2)=-x_1-2$.
It can be checked that all these functions  are convex and the constrained set is nonempty.
The true optimal solution and optimal value for this problem are $(x_1^*, x_2^*)=(1, 2)$ and $\tilde f(x_1^*, x_2^*)=172.41$ respectively.

\begin{figure}[htbp]
\centering
\subfigure{
\begin{minipage}[b]{0.38\textwidth}
\centering
\includegraphics[width=1\textwidth]{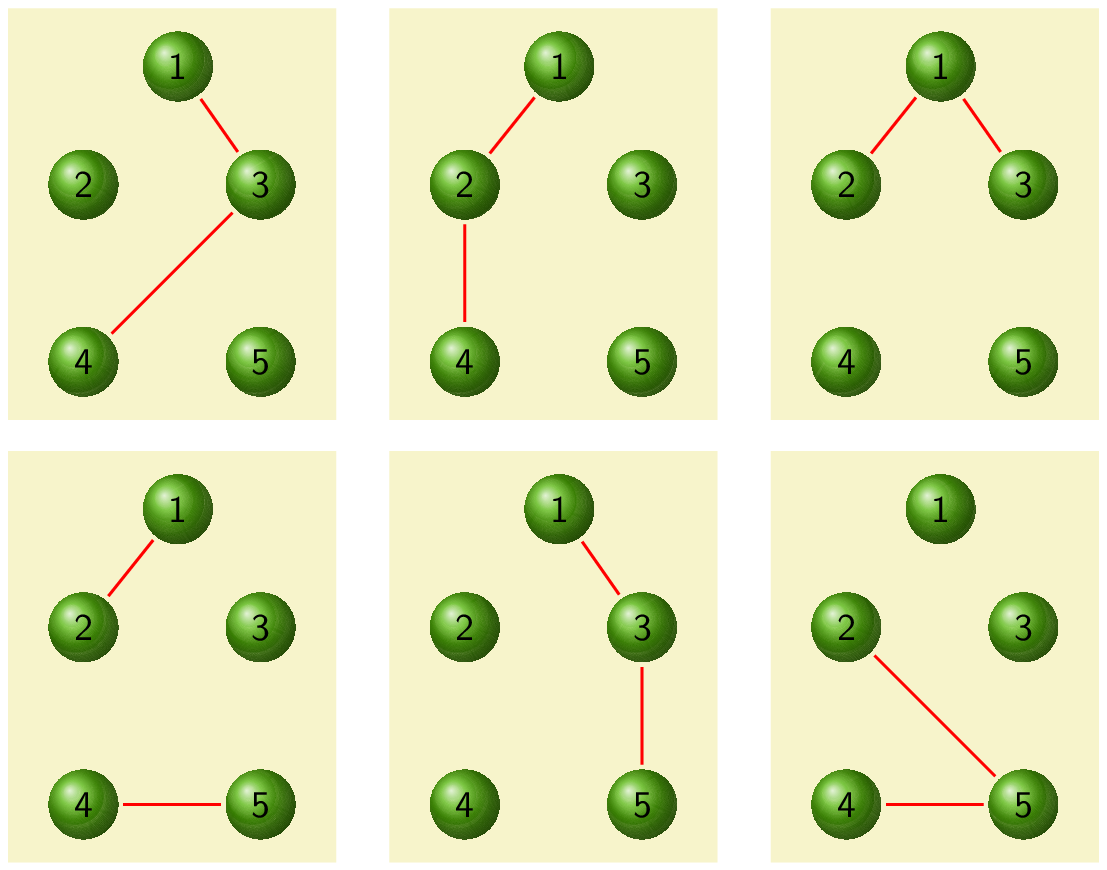}
\end{minipage}%
}%
\subfigure{
\begin{minipage}[b]{0.62\textwidth}
\centering
\includegraphics[width=0.95\textwidth,height=1.5cm]{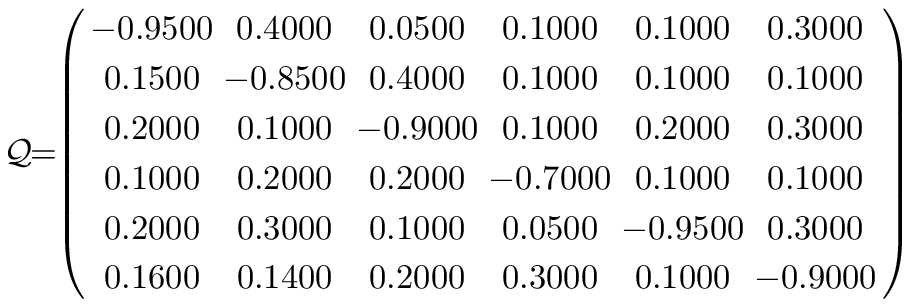}

\includegraphics[width=1\textwidth]{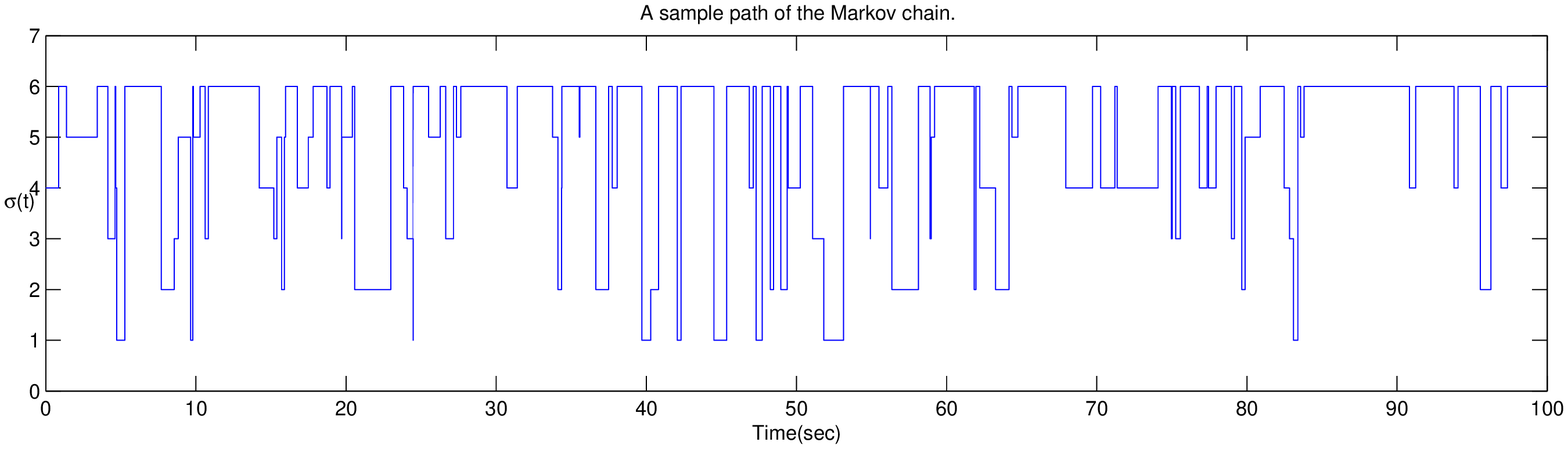}
\centering
\end{minipage}%
}%
\caption{A switching network for the optimization algorithm used in this paper: six possible graphs among which the network switches are shown on the left;  the switching Markov chain is given on the right, where the generator is given on the upper right  and the Markov Switching is plotted on the bottom right.}
\label{FigSixFigure}
\end{figure}

We now use the DCO algorithm \eqref{algorithm} to help check the results.
Let $x_i\in \mathbb{R}^2$ be the state of agent $i\in \{1, \cdots, 5\}$, and its dynamics obeys the algorithm \eqref{algorithm}.
Referring to Figure \ref{FigSixFigure}, the coupling of  five agents   forms a network which is modeled by  a stochastically switching among six possible undirected graphs
$\{\mathcal{G}_1, \cdots, \mathcal{G}_6\}$ (left of Fig. \ref{FigSixFigure}),  with the switching rule described by a continuous-time Markov chain $\sigma: [0, +\infty) \ra \{1, \cdots, 6\}$, whose infinitesimal generator $ \mathcal {Q}$ is given in top right of Figure \ref{FigSixFigure}
which is obviously egordic and the invariant measure can be calculated as
$\pi=(0.1443, 0.2000,$ $0.1882, 0.1652, 0.1132, 0.1891)$, and  whose sample path of $\s$ is shown in bottom right of Figure \ref{FigSixFigure}.
Obviously,  all  graphs are very "sparsely connected",  implying that less communication resources are required at each time. This advantage is more obvious when the number of agents is large. In this sense, switching networks can save communication resources.

\begin{figure}[htbp]
\centering{\includegraphics[width=8.8 cm]{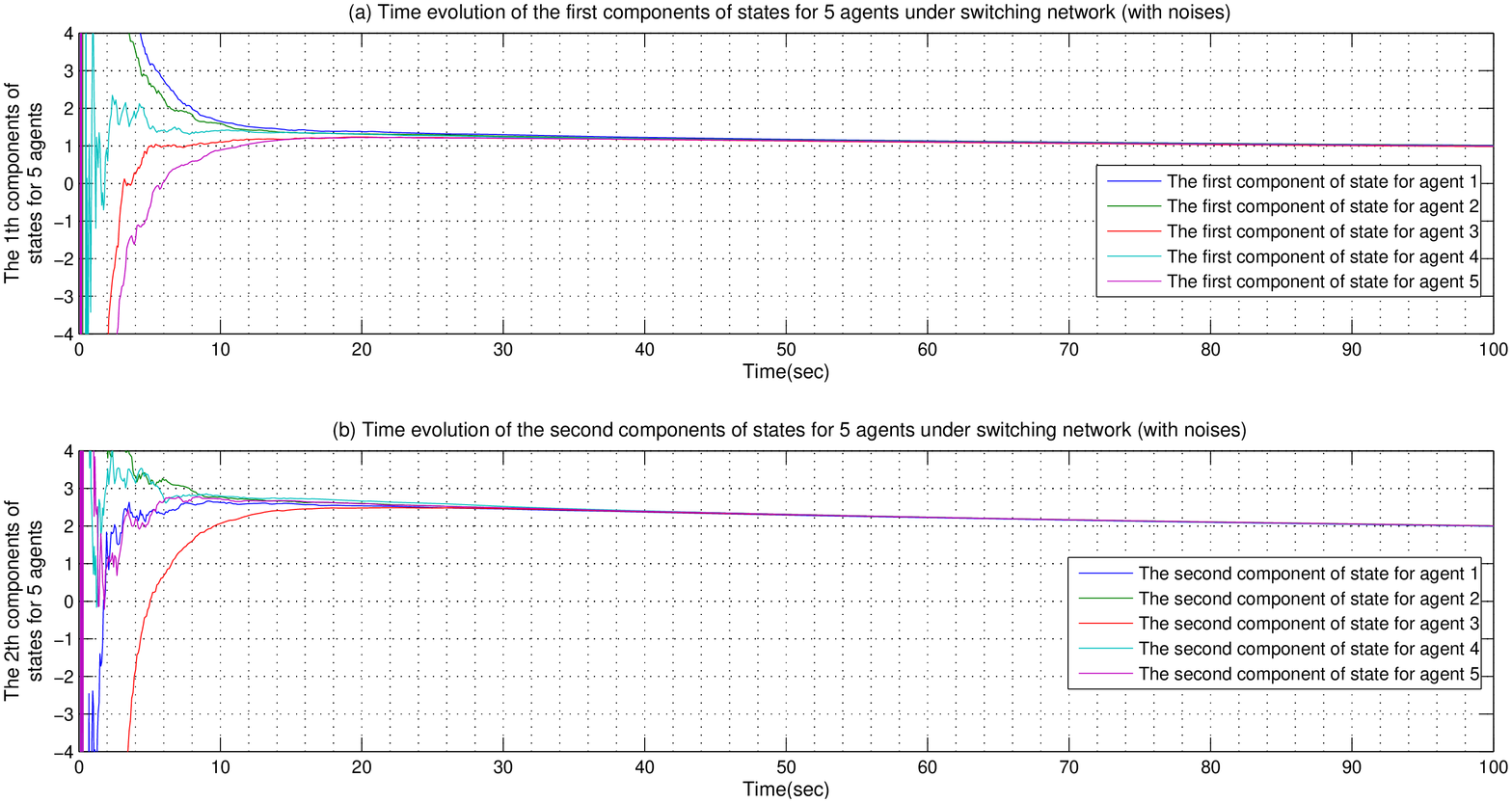}}
\caption{Time evolution of the $x_i$  for 5  agents: (a) The first components  of $x_i's$ converge to $1$; (b) The second components of  $x_i's$ converge to $2$, both almost surely. }
\label{OS}
\end{figure}

For simulation, we chose the noise intensities  $\sigma_{ij}=1$,
the parameters $\eta_{ij}=1$, the coupling strength $c=2$,
and the initial states of five agents  as $x_1(0)\!=\!(-2, 4)^T$, $x_2(0)\!=\!(-3, 3)^T$, $x_3(0)\!=\!(1, -2)^T$, $x_4(0)\!=\!(4, 2)^T$, $x_5(0)\!=\!(-3, -4)^T$, $\theta_1(0)=\cdots=\theta_5(0)=1$, $\lambda_1(0)\!=\!3$, $\lambda_2(0)\!=\!3$, $\nu(0)\!=\!3$. The time evolution of the $\hx$-states for five agents are illustrated in Figure \ref{OS},  where the first component of each state  $x_i$ asymptotically converges to $1$ almost surely (subfigure (a)) and the second component of each state $x_i$ asymptotically converges to $2$ almost surely (subfigure (b)). Therefore, each state $x_i$ of the 5 agents converges to the optimal solution $(1, 2)$ almost surely. Due to space limitation, the time evolutions of the states for $\htheta, \hlambda, \hnu$ are not plotted here.


\section{Conclusion}

This paper has proposed a novel multi-scale method for the  distributed
convex optimization problem with constraints, and presented in a unified
 framework to address challenging issues like optimization constraints,
  communication noises and stochastic networks.
Rigorous convergence analysis is given for the proposed algorithm
 thanks to the use of  Lyapunov arguments, Kolmogorov backward
equation, and Ito formula.
To overcome the technical obstacle in computing the
projection in the presence of optimization constraints, we design a projection-free, smooth optimization dynamics
 for easier analysis
and simulation. As a consequence, a major advantage of the proposed
multi-scale method presented in this paper is that it generalizes previous
distributed convex optimization algorithms from a fixed
network topology to the switching case.
Also, the stochastic averaging in this paper is a generalization of the
deterministic averaging in
our earlier works, again thanks to the multi-scale method
used in this paper.

\appendix

\section*{Appendix A:  Derivation of The SDE  \eqref{X}} \label{appendix1}


To characterize the noisy term $\omega_i:=\sum\nolimits_{j\in \mathscr{N}_i(t)}\sigma_{ji}\xi_{ji}(x_j-x_i)$ in equations \eqref{algorithma}-\eqref{add1},  define $\xi_i=[\xi_{1i}, \xi_{2i}, \cdots, \xi_{Ni}]^T$, $i=1, \cdots, N$. Also, for each switching mode $s\in \mathbb S$, define $\mathcal M^i_{s}\!=\![a^{i1}_{s}\sigma_{1i}(x_1\!-\!x_i), \cdots, a^{iN}_{s}\sigma_{Ni}(x_N-x_i)]$.
Then  the noisy term $\omega_i$  above  can be written as $\omega_i=\mathcal M^i_{s} \xi_i$. Therefore, the equation
\eqref{algorithma} becomes
\begin{align*}
 d x_i\!=\!\!\left[c\!\sum \nolimits_{\scalebox{0.7}[0.8]{$j\!\in\!\! \mathscr{N}_i(t)$}}(x_j\!-\!x_i)
      \!-\!\theta_i\!-\!\nabla f_i(x_i)\!-\!\!\sum \nolimits_{j}^{r_i}\!\!\lambda_{ij}\!\nabla g_{ij}(x_i) \!-\!\sum \nolimits_{j}^{s_i}\!\!\nu_{ij}\nabla\! g_{ij}(x_i)\right]\!dt
   \!+\!c \mathcal M^i_{{\sigma(t)}} dw_i,
\end{align*}
where $\xi_i dt= d w_i$ with $w_i$  an $N$-dimensional standard Browian motion on the probability space $(\Omega, \mathcal{F}, \mathbb{P}, \{\mathcal{F}_t\}_{t\geq 0})$.
Let $\hx={\rm col}\{x_1, \cdots, x_N\}$. Define the stacked functions $G(\hx)={\rm col}\{{\rm col}\{g_{ij}(x_i)\}_{j=1}^{r_i}\}_{i=1}^N \in \mb{R}^r$, $H(\hx)={\rm col}\{{\rm col}\{h_{ij}(x_i)\}_{j=1}^{s_i}\}_{i=1}^N \in \mb{R}^s$
and the stacked gradients $\nabla F(\hx)={\rm col}\{\nabla f_{i}(x_i)\}_{i=1}^{N}\in \mb{R}^{nN}$, $\nabla G(\hx)={\rm col}\big\{{\rm col}\{\nabla g_{ij}(x_i)\}_{j=1}^{r_i}\big\}_{i=1}^{N} \in \mb{R}^{rn}$, $\nabla H(\hx)={\rm col}\big\{{\rm col}\{\nabla h_{ij}(x_i)\}_{j=1}^{s_i}\big\}_{i=1}^{N}\in \mb{R}^{sn}$.
For each switching mode $s\!\in\! \mathbb{S}$,   set $\mathcal{M}_s\!=\!{\rm diag}\{\mathcal{M}_s^i\}_{i=1}^{N}$ $\!\in\! \mb{R}^{nN\times N^2}$.
In addition, let  $\hw\!=\!{\rm col} \{w_i\}_{i=1}^{N}  \in \mb{R}^{N^2}$ and define $r\!=\!\sum_{i=1}^{N}r_i$, $s\!=\!\sum_{i=1}^{N}s_i$.
Set
\scalebox{0.9}[1]{$\hlambda\!=\!{\rm col}\{\lambda_1, \cdots, \lambda_N\}\!\in\! \mb{R}^r$} with \scalebox{0.9}[1]{$\lambda_i \!=\!{\rm col}\{\lambda_{i1}, \cdots, \lambda_{ir_i}\}\!\in\! \mathbb{R}^{r_i}$}, and \scalebox{0.9}[1]{$\hnu\!=\!{\rm col}\{\nu_1, \cdots, \nu_N\}\!\in\!\mb{R}^s$} with \scalebox{0.9}[1]{$\nu_i \!=\!{\rm col}\{\nu_{i1}, \cdots, \nu_{is_i}\}\!\in\! \mathbb{R}^{s_i}$}.
Define $\htheta={\rm col} \{\theta_1, \cdots, \theta_N\}\in \mathbb{R}^{nN}$.
With these, the equation \eqref{algorithm} can be  written in a compact form as in \eqref{X}.


\section*{Appendix B: Proof of The Inequality \eqref{e:appendix}} \label{appendix3}

Firstly, for $V_1$, defining $\hh= \nabla F(\hx)+\boldsymbol{\lambda} \odot \nabla G(\hx)+\boldsymbol{\nu} \odot \nabla H(\hx)$ and $\hh^*= \nabla F(\hx^*)+\!\boldsymbol{\lambda}^* \odot \nabla G(\hx^*)+\boldsymbol{\nu}^* \odot \nabla H(\hx^*)$ and noting $\htheta^*=-\hh^*$ in view of \eqref{conditionTheta}, the action of the infinitesimal operator $\mc A$  on   $V_1$ can be calculated as

\begin{align}\label{AV1}
  \mc{A}V_1
  &=\underbrace{(\hx\!-\!\hx^*)^T[-c\boldsymbol{\mathcal{L}}\hx-\htheta-\n F(\hx)]+\frac{1}{2}c^2tr(\mathcal{M}^T\mathcal{M})}_{A}-\sum_{i=1}^{N}\sum_{j=1}^{r_i}\lambda_{ij}(x_i\!-\!x^*)^T \nabla g_{ij}(x_i)\nonumber\\
  &\hspace{0.5cm}+\!\underbrace{[(\hx\!-\!\hx^*)\!+\!(\htheta\!-\!\htheta^*)]^T[-(\htheta\!-\!\htheta^*)\!-\!(\hh\!-\!\hh^*)]}_{B}-\!\sum_{i=1}^{N}\!\sum_{j=1}^{s_i}\!\nu_{ij}(x_i-x^*)^T \nabla h_{ij}(x_i)
\end{align}
The trace ${\rm tr}(\mathcal{M}^T\mathcal{M})$  brings difficulty to the convergence analysis. To get rid of this difficulty, we  give an estimation of ${\rm tr}(\mathcal{M}^T\mathcal{M})$ as follows,
\begin{align}\label{MM}
\hspace{-0.1cm}{\rm tr}(\mathcal{M}^T\!\mathcal{M})
&=\sum \nolimits_{i=1}^N {\rm tr}[(\mathcal M^i)^T\mathcal M^i]=\sum \nolimits_{i=1}^N \sum\nolimits_{j=1}^N(a^{ij})^2(\sigma_{ji})^2(x_j-x_i)^T(x_j-x_i)\nonumber\\
&\leq \kappa^2\sum\nolimits_{i=1}^N \sum\nolimits_{j=1}^Na^{ij}(x_j-x_i)^T(x_j-x_i) \!=\!\kappa^2 \hx^T \!\boldsymbol{\mathcal{L}} \hx.
\end{align}
As for Part A in \eqref{AV1} , noting $\hx^T \!\boldsymbol{\mathcal{L}} \hx\!=\! (\hx\!-\!\hx^*)^T \boldsymbol{\mathcal{L}} (\hx\!-\!\hx^*)$ and $\hbar\!=\!\scalebox{1}[0.9]{$\frac{1}{2}$}(c\!-\!\frac{1}{2}c^2\kappa^2)$ (c.f. Eq. \eqref{psi}), the part A can be estimated as $A\!\leq\! (\hx\!-\!\hx^*)^T\![-\hbar\boldsymbol{\mathcal{L}}\hx\!-\!\htheta\!-\!\n F(\hx)]-\hbar(\hx\!-\!\hx^*)^T \boldsymbol{\mathcal{L}}(\hx\!-\!\hx^*)$.
Noting that the square bracket  is exactly the minus gradient of $\Psi(\hx, \htheta)$  (c.f. Eq. \eqref{psi}) with respect to $\hx$ and denoting by $\lambda^{\mathcal G}_2$ the smallest nonzero eigenvalue of the Laplacian $\mathcal L$, one obtains $A \!\leq\! (\hx^*\!-\!\hx)^T \n_{\!\hx}  \Psi(\hx, \htheta)\!-\!\hbar \lambda^{\mathcal G}_2(\hx\!-\!\hx^*)^T(\hx\!-\!\hx^*)\!\leq\! \Psi(\hx^*, \htheta)\!-\!\Psi(\hx,\htheta)\!-\!\hbar \lambda^{\mathcal G}_2(\hx\!-\!\hx^*)^T(\hx\!-\!\hx^*)$, where the last inequality uses the fact that $ \Psi(\hx, \htheta)$ is convex in its first argument.
As for part  B in \eqref{AV1} , we can prove the  inequality $B\leq (\hx-\hx^*)^T(\htheta-\htheta^*)+(\hx-\hx^*)^T(\hx-\hx^*)$
by rewriting it into a quadratic form in terms of $\hx-\hx^*, \htheta-\htheta^*, \hh-\hh^*$ and by showing the corresponding matrix to be semi-positive definite.
In view of  $(\hx-\hx^*)^T(\htheta-\htheta^*)=\Psi(\hx, \htheta)-\Psi(\hx, \htheta^*)$,  one obtains
$A+B\leq \left[ \Psi(\hx^*, \htheta)-\Psi(\hx, \htheta^*)\right]-(\hbar \lambda^{\mathcal G}_2-1)(\hx-\hx^*)^T(\hx-\hx^*)$.
Now, $\mathcal{A}V_1$ can be calculated as
\begin{align*}
 \mathcal{A}V_1
&\leq   \Psi(\hx^*, \htheta)- \Psi(\hx, \htheta^*)\!-\!(\hbar \lambda^{\mathcal G}_2-1)(\hx\!-\!\hx^*)^T(\hx\!-\!\hx^*)+\!\sum\nolimits_{i=1}^{N}\sum\nolimits_{j=1}^{r_i}\lambda_{ij}g_{ij}(x^*)\nonumber\\
&\!-\!\sum\nolimits_{i=1}^{N}\sum\nolimits_{j=1}^{r_i}\lambda_{ij}g_{ij}(x_i)+\sum\nolimits_{i=1}^{N}\sum\nolimits_{j=1}^{s_i}\nu_{ij}h_{ij}(x^*)\!-\!\sum\nolimits_{i=1}^{N}\sum\nolimits_{j=1}^{s_i}\nu_{ij}h_{ij}(x_i),
\end{align*}
where we use the convexity of the functions $g_{ij}$ and $h_{ij}$.

Secondly, we  calculate the action of the infinitesimal operator $\mc A$  on   $V_2+V_3$.
Firstly note that  the function $V_3$ can be calculated by the definition of Bregman divergence (c.f. \cite{bregman1967}) as
 $V_3=\sum\nolimits_{i=1}^{N} \sum\nolimits_{j=1}^{r_i}(\lambda_{ij}-\lambda_{ij}^*)-\sum_{(i,j)\in \Omega}\lambda_{ij}^*(\ln \lambda_{ij}-\ln \lambda_{ij}^*)$.  Therefore,
\begin{align*}
\mathcal{A}V_2 + \mc{A}V_3
&=\sum\nolimits_{i=1}^{N}\sum\nolimits_{j=1}^{r_i} \frac{\eta_{ij}\lambda_{ij}(\lambda_{ij}-\lambda_{ij}^*)}{1+\eta_{ij}\lambda_{ij}} g_{ij}(x_i)+\sum\nolimits_{i=1}^{N}\sum\nolimits_{j=1}^{r_i}  \frac{\lambda_{ij}}{1+\eta_{ij}\lambda_{ij}} g_{ij}(x_i)\\
&\hspace{1.5cm}-\sum\nolimits_{(i,j)\in \Omega} \frac{\lambda_{ij}^*}{1+\eta_{ij}\lambda_{ij}} g_{ij}(x_i)=\sum\nolimits_{i=1}^{N}\sum\nolimits_{j=1}^{r_i} (\lambda_{ij}-\lambda_{ij}^*)g_{ij}(x_i).
\end{align*}
Furthermore, the action of the infinitesimal operator $\mc A$  on   $V_4$ can be easily calculated as $\mathcal{A}V_4=\sum\nolimits_{i=1}^{N}\sum\nolimits_{j=1}^{s_i} (v_{ij}-v_{ij}^*) h_{ij}(x_i)$.
Collecting above results for $\mc{A}V_i, i=1,2,3,4$ and recalling the definition of  $ \Phi$ in \eqref{barlagrange},  one obtains the  inequality \eqref{e:appendix}.

\section*{Appendix C: Proof of The Saddle Point Conditions  \eqref{SaddlePointCondition}} \label{appendix4}

Due to convexity and affinity,  the following results hold
\begin{align*}
&\Sigma_{i=1}^N f_i(x_i)+\hbar \hx \hL \hx \geq \Sigma_{i=1}^N f_i(x^*)+\Sigma_{i=1}^N \n f_i(x^*)(x_i-x^*),\\
& g_{ij}(x_i)\geq g_{ij}(x^*)+\n g_{ij}(x^*)(x_i-x^*),\\
& h_{ij}(x_i) = h_{ij}(x^*)+ \n h_{ij}(x^*)(x_i-x^*),\\
& \hx-\hx^*=\hx-\hx^*.
\end{align*}
Multiplying them by $1,  \lambda_{ij}^*, \nu_{ij}^*, \htheta^*$ respectively (for the forth equality we use $(\hx-\hx^*)^T \htheta^*$ for multiplication) and adding gives
$\Phi(\hx, \htheta^*, \hlambda^*, \hnu^*)\!\geq\!\Phi(\hx^*, \htheta^*, \hlambda^*, \hnu^*)\!+\![\htheta^*\!+\! \nabla F(\hx^*) \!+\!\boldsymbol{\lambda}^* \!\odot\! \nabla G(\hx^*)\!+\!\boldsymbol{\nu}^* \!\odot\! \nabla H(\hx^*)]^T(\hx\!-\!\hx^*)$. In view of the equation \eqref{conditionTheta}, the terms in the square bracket sum to be zero. Therefore,  $\Phi(\hx, \htheta^*, \hlambda^*, \hnu^*) \geq  \Phi(\hx^*, \htheta^*, \hlambda^*, \hnu^*)$.

On other hand, noting that
${\color{red}\lambda_{ij}g_{ij}(x_i^*)\leq 0}=\lambda_{ij}^*g_{ij}(x_i^*)$ and $h_{ij}(x^*)=0$, it can be directly checked that $\Phi(\hx^*, \htheta, \hlambda, \hnu)\leq \Phi(\hx^*, \htheta^*, \hlambda^*, \hnu^*)$.
This inequality, together with the inequality derived in last paragraph, gives rise to the saddle point condition  \eqref{SaddlePointCondition}.

\section*{Appendix D: Proof of  ``$\mathcal{A}V=0 \Rightarrow (\hx, \htheta, \hlambda, \hnu)=(\hx^*, \htheta^*, \hlambda^*, \hnu^*)$''} \label{appendix4}

Letting $\mathcal{A}V=0$ gives $(\hx-\hx^*)^T \,(\hx-\hx^*)=0$ and $\Phi(\hx^*, \htheta, \hlambda, \hnu)= \Phi(\hx, \htheta^*, \hlambda^*, \hnu^*)$, where the former implies $\hx=\hx^*$ and the latter, together with the fact that $(\hx^*, \htheta^*, \hlambda^*, \hnu^*)$ is a saddle point of the Lagrangian $ \Phi$,  implies $ \Phi(\hx^*, \htheta, \hlambda, \hnu)\!=\! \Phi(\hx, \htheta^*, \hlambda^*, \hnu^*)\!=\!\Phi(\hx^*, \htheta^*, \hlambda^*, \hnu^*)$.
Recall that $\hx^*=\mathbf 1_N \otimes x^*$ with $x^*$ being the optimal solution which satisfies $g_{ij}(x^*)\leq 0$ in the KKT condition \eqref{KKTa}.  Inserting $x_i=x^*$ into \eqref{algorithmb} yields $\dot\lambda_{ij}=\frac{\lambda_{ij}g_{ij}(x^*)}{1+\eta_{ij}\lambda_{ij}}$.
If $g_{ij}(x^*)=0$, then  $\dot \lambda_{ij}=0$ which implies that  $\lambda_{ij}(t)$ stays positive for all $t\geq 0$  since the initial value of $\lambda_{ij}$ is chosen to be positive.  If $g_{ij}(x^*)<0$, then with $g_{ij}(x^*)=-a<0$, one has   $\dot\lambda_{ij}=\frac{-a\lambda_{ij}}{1+\eta_{ij}\lambda_{ij}}$  whose trajectory can be shown by elementary analysis  as $\lambda_{ij}(t)\geq 0$ for $t\geq 0$  since the initial value of $\lambda_{ij}$ is chosen to be positive.   In short,  in both cases, for the equation \eqref{algorithmb} with $x_i=x^*$, one has that $\lambda_{ij}(t)\geq 0, \forall t\geq 0$ provided the initial value is positive, and consequently $\dot\lambda_{ij}(t)\leq 0$ since $g_{ij}(x^*)\leq 0$. Therefore,  $\lambda_{ij}\geq \lambda_{ij}^*$.
On the other hand, the fact $ \Phi(\hx^*, \htheta, \hlambda, \hnu)\!=\! \Phi(\hx^*, \htheta^*, \hlambda^*, \hnu^*)$ yields $\sum_{i=1}^{N}\sum_{j=1}^{r_i}(\lambda_{ij}\!-\!\lambda_{ij}^*)g_{ij}(x^*)$ $\!=\!0$.
Since $g_{ij}(x^*)\leq 0$ and $\lambda_{ij}\geq \lambda_{ij}^*$, each sum in this summation is non-positive and therefore $(\lambda_{ij}-\lambda_{ij}^*)g_{ij}(x^*)=0$, i.e.,  $\lambda_{ij}g_{ij}(x^*)=\lambda_{ij}^*g_{ij}(x^*)=0$.
This fact implies that the right hand side of equation \eqref{algorithmb} with $x_i=x^*$ is zero and thus $ \lambda_{ij}=\lambda_{ij}^{\dag}$ for some constant $\lambda_{ij}^{\dag}\!\geq\!0$. Thus $\lambda_{ij}^{\dag}g_{ij}(x^*)\!=\!0$.    Inserting $\hx^*\!=\!\mathbf 1_N\otimes x^*$ into equations \eqref{add1} and \eqref{algorithmc} gives that $\theta_i=\theta_i^{\dag}$ and $\nu_{ij}=\nu_{ij}^{\dag}$ for some constants $\theta_i^{\dag}$ and $\nu_{ij}^{\dag}$.  Inserting $(\hx^*, \theta_i^{\dag}, \lambda_{ij}^{\dag}, \nu_{ij}^{\dag})$ into \eqref{algorithma} gives $\sum_{i=1}^N\nabla f_i(x^*)\!+\!\sum_{i=1}^N\sum_{j}^{r_i}\lambda_{ij}^{\dag}\nabla g_{ij}(x^*)\!+\!\sum_{i=1}^N\sum_{j}^{s_i}\nu_{ij}^{\dag}\nabla h_{ij}(x^*)\!=\!0$.
In conclusion, the KKT conditions \eqref{KKT} are satisfied at the point $(\hx^*, \hlambda^\dag, \hnu^\dag)$.
Then  by uniqueness of multipliers in \eqref{uniqueMultiplier},  $\lambda_{ij}^{\dag}=\lambda_{ij}^*$, $\nu_{ij}^{\dag}=\nu_{ij}^*$.

Also, by differentiating  both sides of $ \Phi(\hx, \htheta^*, \hlambda^*, \hnu^*)\!=\! \Phi(\hx^*, \htheta^*, \hlambda^*, \hnu^*)$ with respect to $\hx$ and then by enforcing $\hx=\hx^*$, one obtains $\htheta\!+\! \nabla F(\hx^*) \!+\!\boldsymbol{\lambda}^* \!\odot\! \nabla G(\hx^*)\!+\!\boldsymbol{\nu}^*\! \odot\! \nabla H(\hx^*)=0$. This, combined with  equation \eqref{conditionTheta}, gives $\htheta=\htheta^*$.
In conclusion, letting $\mathcal{A}V=0$ gives rise to $(\hx, \htheta, \hlambda, \hnu)=(\hx^*, \htheta^*, \hlambda^*, \hnu^*)$.



\begin{acknowledgements}
The first author would like to thank Professor Zhong-Ping Jiang for his useful discussions on this paper when the first author visited New York University. This work is supported by   the NNSF of China under the grants 61663026, 61473098, 61563033, 11361043,61603175.
\end{acknowledgements}

\appendix  

\bibliography{ref}


%
%
%
%
%
%
%

\end{document}